\documentclass[11pt,a4paper]{article}
\usepackage[T1]{fontenc}
\usepackage[utf8]{inputenc}
\usepackage{lmodern,amsmath,amssymb,amsthm,booktabs,graphicx,tikz,microtype}
\usepackage[margin=25mm]{geometry}
\usepackage[font=small,labelfont=bf]{caption}
\usepackage[hidelinks]{hyperref}

\newtheorem{proposition}{Proposition}

\begin{document}

\begin{center}
{\bf \Large A 4-regular graph with \boldmath{4/5} of its vertices\\ satisfying the \v{S}olt\'es property\par}
\end{center}

\begin{center}
A.~A.~Dobrynin

\smallskip

\textit{Sobolev Institute of Mathematics, Siberian Branch of the \\ Russian Academy of Sciences, Novosibirsk, 630090, Russia}

email: dobr@math.nsc.ru
\end{center}

\begin{abstract}
The Wiener index $W(G)$ of a connected graph $G$ is the sum of distances
between all unordered pairs of its vertices. A vertex $v$ is called a
\v{S}olt\'es vertex if $G-v$ is connected and $W(G-v)=W(G)$. Since the cycle
$C_{11}$ is the only known graph in which every vertex is a \v{S}olt\'es vertex,
a natural problem is to construct other graphs with as large a proportion
of such vertices as possible. We construct a $4$-regular graph with
\v{S}olt\'es fraction $4/5$, improving the previous value of $2/3$.

\smallskip


\smallskip

\end{abstract}

\section*{Introduction}

All graphs considered are simple, undirected and connected.
The vertex set of a graph $G$ is denoted by $V(G)$, and its \emph{order} is
$|V(G)|$. The \emph{girth} of $G$ is the length of its shortest cycle. 
The distance $d_G(u,v)$ between two vertices $u$ and $v$ is the
length of a shortest path joining them. 
The eccentricity of a vertex $v$ in a graph $G$ is
 the maximum distance between $v$ and any other vertex of $G$.
The \emph{radius} and \emph{diameter} of a graph $G$ are the 
minimum and the maximum vertex eccentricity, respectively.
The \emph{transmission} of a
vertex $v$ of $G$ is the sum of the distances from $v$ to all other
vertices of $G$, that is, $T_G(v)=\sum_{u\in V(G)}d_G(v,u)$. The \emph{Wiener index}
of $G$ is the sum of distances over all unordered pairs of vertices,
$$
W(G)=\sum_{\{u,v\}\subseteq V(G)}d_G(u,v)
=\frac12\sum_{v\in V(G)}T_G(v).
$$

The Wiener index was introduced as a structural descriptor for
molecular graphs~\cite{Wiener}. Its mathematical properties and 
numerous applications, in particular in organic chemistry,
have been extensively studied; see, for example,
\cite{DobryninEntringerGutman,DobryninGutmanKlavzarZigert,
KnorSkrekovskiTepeh2016,KnorSkrekovskiTepeh2024,
Trinajstic1992,WagnerWang2018}.

 A vertex $v\in V(G)$ is called a \emph{\v{S}olt\'{e}s vertex} if $G-v$ is connected and $W(G-v)=W(G)$.
In 1991, \v{S}olt\'es~\cite{Soltes} posed the problem of finding graphs
in which all vertices are  \emph{\v{S}olt\'es vertices}.
He also gave the cycle $C_{11}$ as the first example.
Since $C_{11}$ is the only known example, various relaxations of this
problem have been studied, including
the question of how large a proportion of the vertices can have
this property. 
The \emph{\v{S}olt\'es fraction} $R_0(G)$  is the proportion of  vertices of $G$ that are \v{S}olt\'es vertices.
A number  of graph families with prescribed numbers of \v{S}olt\'es vertices and
graphs from restricted classes were constructed in 
\cite{BasicKnorSkrekovski,BokJedlickovaMaxova2019,BokJedlickovaMaxova2021,
Dobrynin2023,DobryninVorobev,HuZhuWuShaoFahad2021,
KnorMajstorovicSkrekovskiAMC,KnorMajstorovicSkrekovskiDAM}.
The \v{S}olt\'es fraction of graphs from these studies does not exceed 2/3.
Akhmejanova et al.~\cite{AkhmejanovaEtAl} presented the first graph of order $12$
with eight \v{S}olt\'es vertices and  the fraction $2/3$. Cambie~\cite{Cambie2026} described eight graphs with
\v{S}olt\'es fraction $2/3$, including this example. He also proposed a conditional
construction of graphs with \v{S}olt\'es fraction approximately $5/7$.

In this note we construct a $4$-regular graph of order $2400$ with
exactly $1920$ \v{S}olt\'es vertices, and hence with \v{S}olt\'es fraction $4/5$.

\section{The base graph}\label{sec:base}

To construct a 4-regular graph with \v{S}olt\'es fraction $4/5$, we use the
base graph $H$ described below, which will be transformed as explained
in the next section.

Let $H$ denote the graph \texttt{AT4val[480,150]} in the census of
connected tetravalent arc-transitive graphs of order at most $640$
compiled by Poto\v{c}nik, Spiga and Verret~\cite{PSV2013,PSV2015,census}.
The graph $H$ is connected, bipartite and $4$-regular, with $480$
vertices and $960$ edges. Its girth, diameter and radius are all $8$.
The graph $H$ is a Cayley graph of the group identified as 
SmallGroup(480,955) in the GAP Small Groups Library~\cite{SmallGrp}.
The automorphism group of $H$ has order $1920$ and
acts regularly on the arcs. In particular, $H$ is vertex-transitive
and edge-transitive.
A drawing of $H$ is shown in Figure~\ref{fig:base} (see the Appendix for other drawings).
The distance partition with respect to a vertex has layers of sizes
$1$, $4$, $12$, $36$, $93$, $154$, $129$, $46$, and $5$. Hence every
vertex has transmission $2414$. For every
$v\in V(H)$, we have $W(H)-W(H-v)=579\,360-577\,598=1\,762>0$.

There are $1080$ cycles of length $8$ in $H$, and each vertex belongs
to $18$ of them. A computational analysis establishes the following property of $H$.

\textbf{Cyclic property.}  {\em For every vertex $v$, four of the six paths of length
$2$ with middle vertex $v$ are contained in exactly one cycle of length
$8$, and the other two, which have no edge in common, are contained in exactly seven such cycles.}

\begin{figure}[h]
\centering
\includegraphics[trim=6.26bp 6.05bp 6.05bp 6.26bp,clip,width=0.93\textwidth,height=0.93\textheight,keepaspectratio]{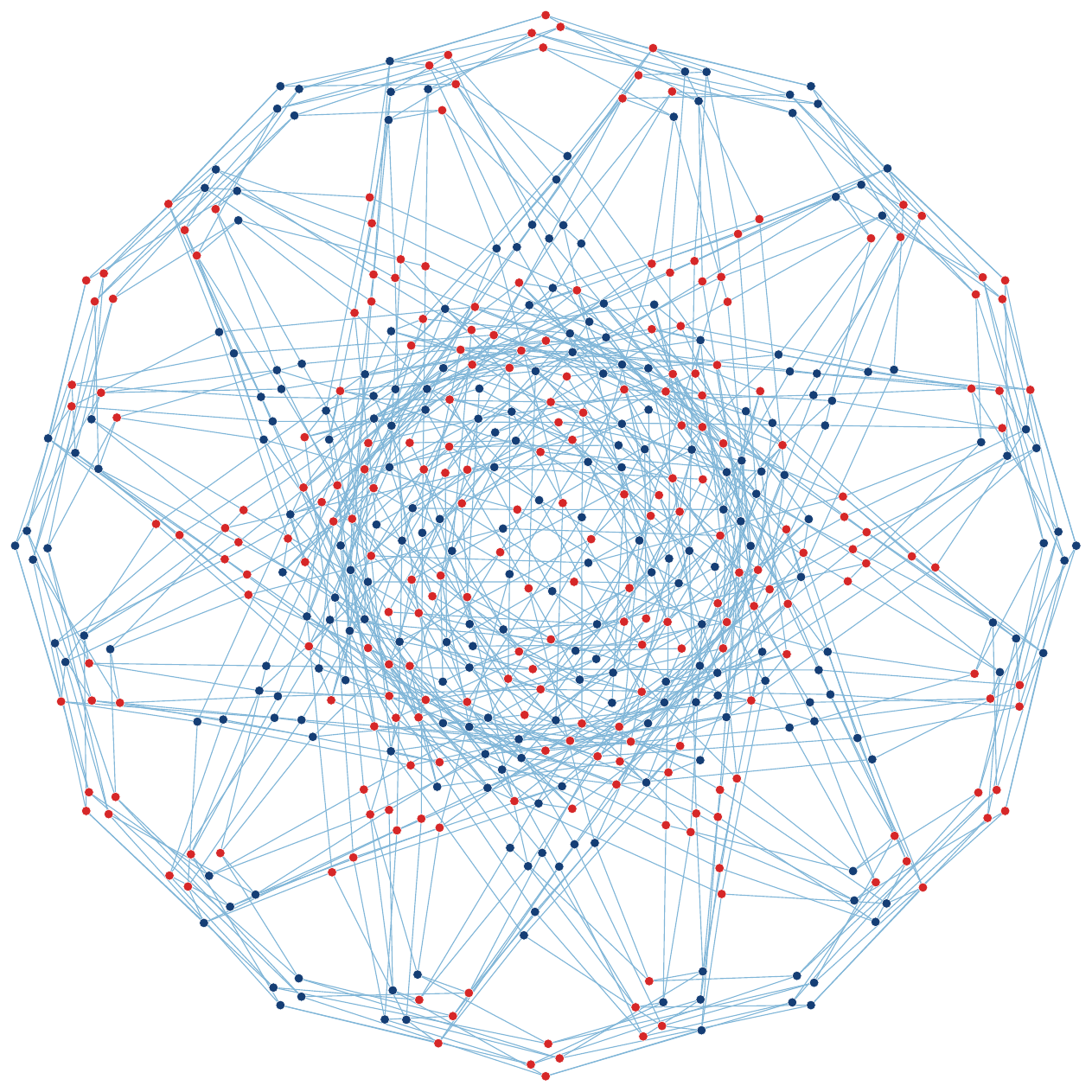}
\caption{The base graph $H$. Red and dark blue indicate its two bipartition classes.}
\label{fig:base}
\end{figure}

\clearpage

\section{Construction of the graph}\label{sec:construction}

Let $H$ be the base graph described above. We construct $G$ from $H$ as follows.

\begin{figure}[t!]
\centering
\includegraphics[width=\textwidth,height=\textheight,keepaspectratio]{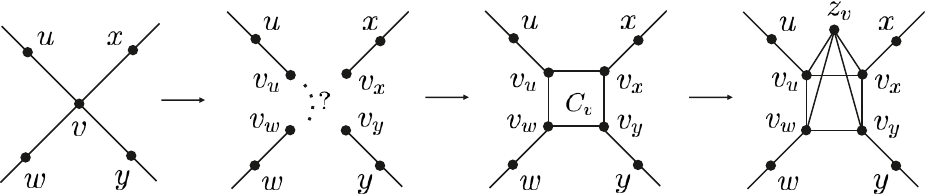}
\caption{Transformation of the base graph $H$ into $G$.}
\label{fig:construction}
\end{figure}

\begin{enumerate}
\item
Replace each vertex $v$ of $H$ by four vertices $v_u$, $v_w$, $v_x$, and $v_y$,
where $u$, $w$, $x$, and $y$ are the neighbors of $v$,
as shown in Figure~\ref{fig:construction}.
Each original edge $uv$ becomes the edge $u_v v_u$.
These new vertices have degree $1$ and correspond to the arcs of $H$.

\item
For each vertex $v$ and each pair of distinct neighbors $u$ and $w$ of $v$, join $v_u$ and $v_w$ 
if and only if the path $u v w$ lies in exactly one cycle of length 8 in $H$. 
By the cyclic property, for every $v$ exactly two of the six pairs are excluded by this rule, 
and these two pairs are disjoint. Hence vertices $v_u, v_w, v_x, v_y$  induce a 4-cycle $C_v$.

\item
For each vertex $v$, add an auxiliary vertex $z_v$ and join it to all four vertices of $C_v$.
\end{enumerate}

The graph $G$ is connected and $4$-regular, with order
$|V(G)|=5\cdot480=2400$ and size $|E(G)|=960+8\cdot480=4800$.
Its diameter and radius are both $18$, and the girth is $3$.

Every automorphism of $H$ preserves the rule that defines the added edges,
and therefore extends to an automorphism of $G$. The vertices of the
cycles $C_v$ correspond to the arcs of $H$, and the auxiliary
vertices correspond to the vertices of $H$. Since $H$ is arc-transitive,
the group induced by $\operatorname{Aut}(H)$ acts transitively on each of these two sets.
A vertex of a cycle $C_v$ lies in exactly two triangles, whereas an
auxiliary vertex lies in exactly four. 
Hence no automorphism of $G$ maps a vertex of one set to a vertex of the other,
so these two sets
are precisely the two vertex orbits of $\operatorname{Aut}(G)$, of sizes $1920$
and $480$.

\begin{proposition}\label{prop:four-fifths}
The graph $G$ of order $2400$ has exactly $1920$ \v{S}olt\'es vertices,
and hence its \v{S}olt\'es fraction equals $4/5$.
\end{proposition}

The graph $G-v$ is connected for every $v\in V(G)$. The differences
$\Delta_G(v)=W(G)-W(G-v)$ were computed by the author for all $2400$
vertices using a program written in Delphi 12. The results for the two vertex orbits are shown in
Table~\ref{tab:deletions}.

\begin{table}[htbp]
\centering
\caption{Transmissions and vertex-deletion values for $G$.}
\label{tab:deletions}
\begin{tabular}{@{}lrrrr@{}}
\toprule
Vertex type & Number & $T_G(v)$ & $W(G-v)$ & $\Delta_G(v)$\\
\midrule
Quadrilateral & $1920$ & $27\,289$ & $32\,823\,120$ & $0$\\
Auxiliary & $480$ & $27\,607$ & $32\,795\,513$ & $27\,607$\\
\bottomrule
\end{tabular}
\end{table}

Thus $W(G)=(1920\cdot27\,289+480\cdot27\,607)/2=32\,823\,120$,
and the \v{S}olt\'es vertices are precisely the $1920$ vertices of the
quadrilaterals. Therefore, $R_0(G)=4/5$.

The dataset containing the edge list of $G$ and its \textit{graph6} and \textit{sparse6}
encodings is available from the author upon request.

\section*{Acknowledgments}

This work was supported by the state contract of the Sobolev Institute
of Mathematics (project number FWNF-2022-0017).

\section*{Declaration of generative AI use}

ChatGPT (OpenAI) was used to assist with mathematical discussions and computational checks. 
All drawings of the base graph $H$ were produced using ChatGPT. 
The author reviewed and revised the material and takes full responsibility for the final content.

\section*{Appendix}

The graph $H$ admits Cayley representations over three pairwise non-isomorphic groups, identified in the Small Groups Library of GAP 
as SmallGroup(480,219), SmallGroup(480,955), and SmallGroup(480,959). The figures below show drawings of $H$ with different rotational symmetries.

\

\

\begin{figure}[h]
\centering
\includegraphics[trim=26.26bp 26.05bp 26.05bp 26.26bp,clip,width=0.96\textwidth,height=0.96\textheight,keepaspectratio]{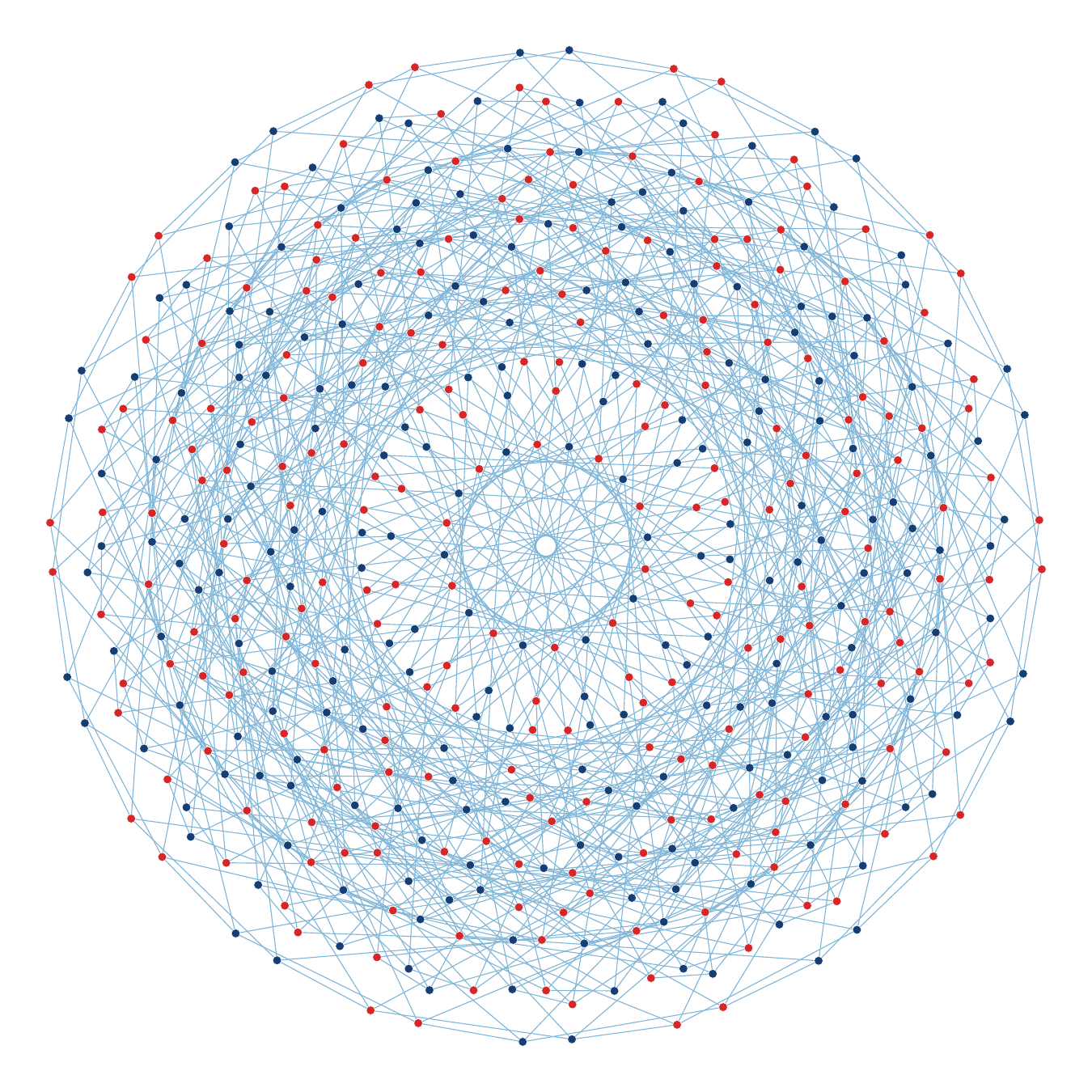}
\caption{An alternative drawing of the base graph $H$, with the same vertex colors.}
\end{figure}

\clearpage

\

\

\

\begin{figure}[h]
\centering
\includegraphics[trim=6.26bp 6.05bp 6.05bp 6.26bp,clip,width=\textwidth,height=\textheight,keepaspectratio]{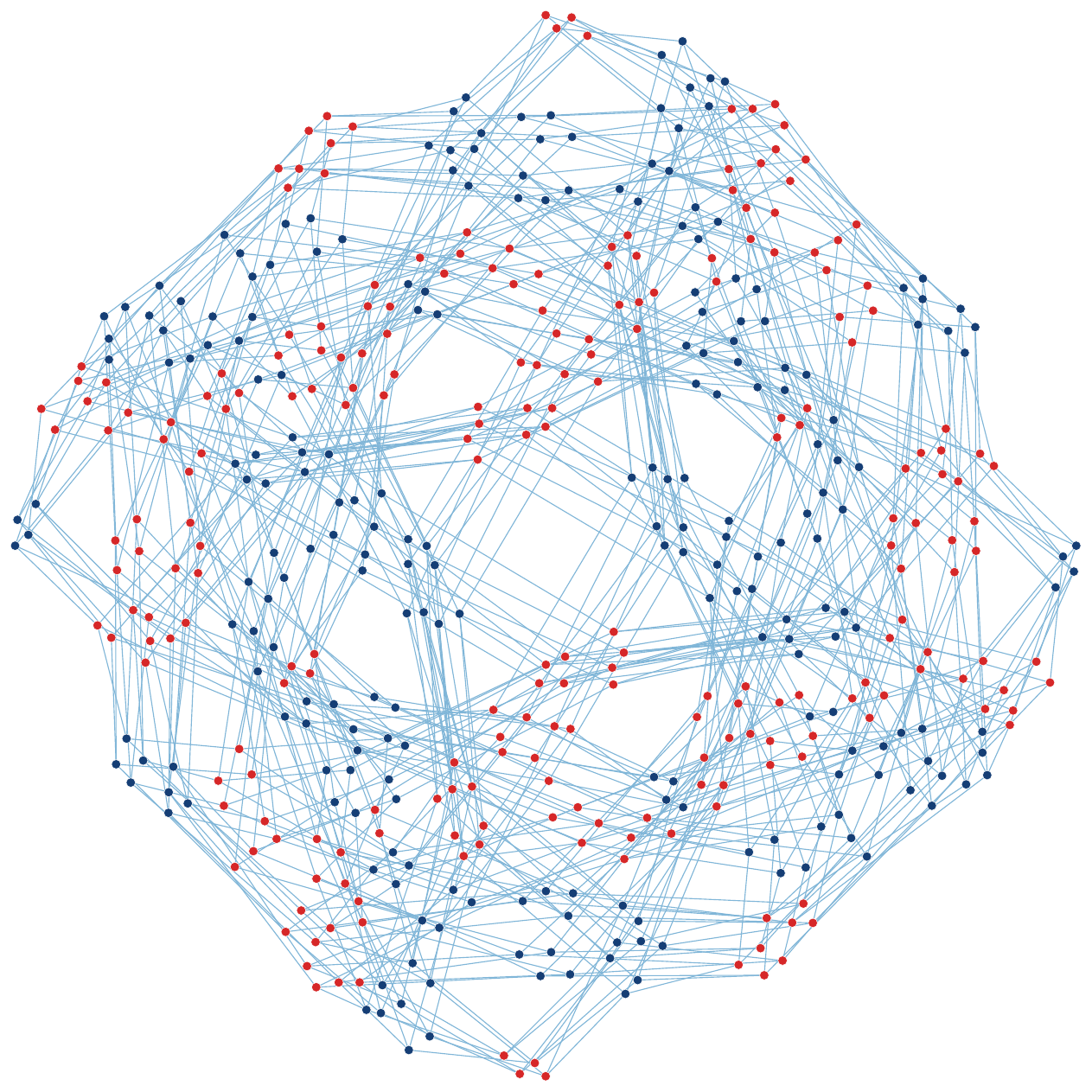}
\caption{An alternative drawing of the base graph $H$, with the same vertex colors.}
\end{figure}

\

\

\clearpage

\

\

\begin{figure}[h]
\centering
\includegraphics[trim=06.26bp 06.05bp 06.05bp 06.26bp,clip,width=\textwidth,height=\textheight,keepaspectratio]{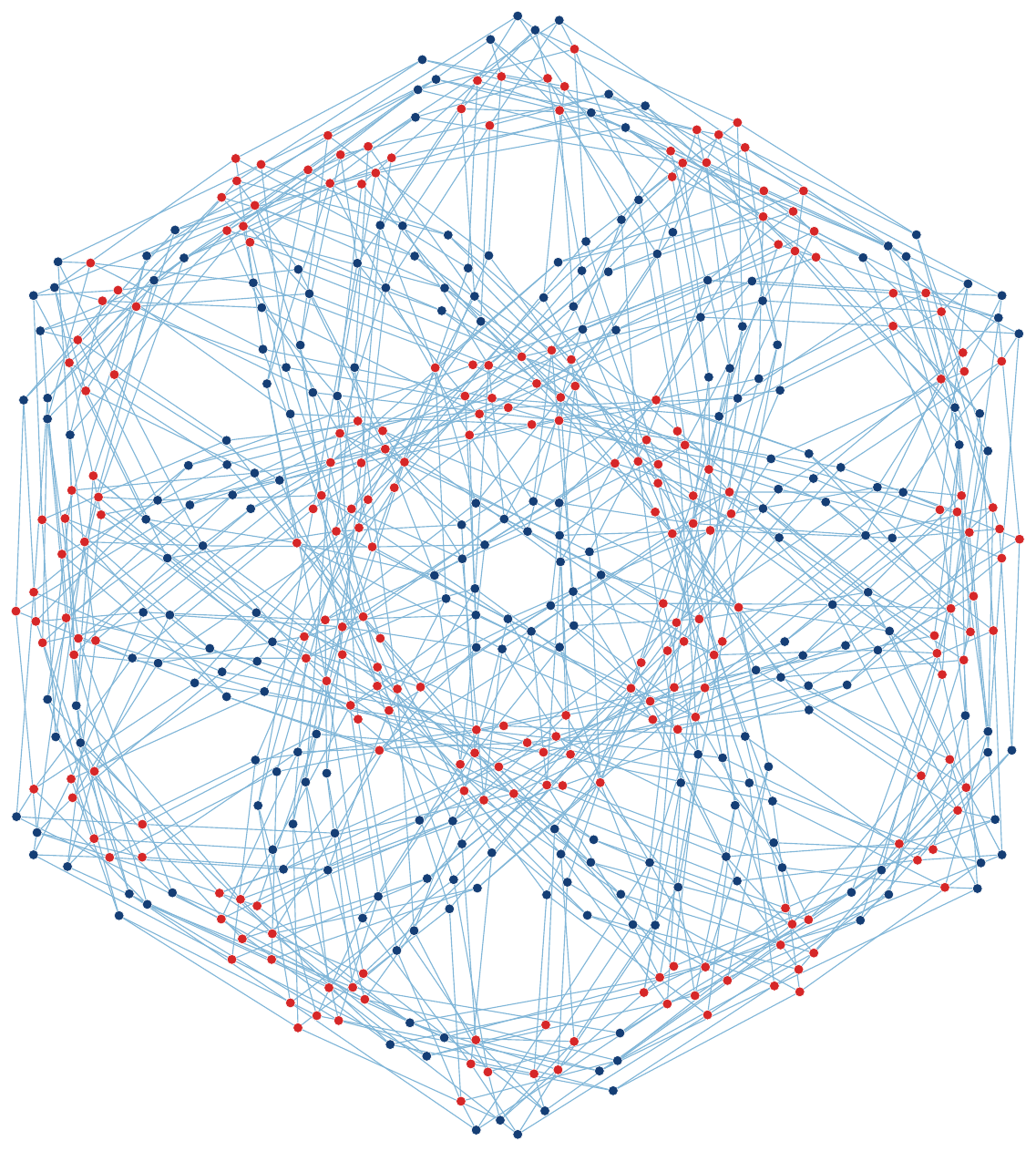}
\caption{An alternative drawing of the base graph $H$, with the same vertex colors.}
\end{figure}

\clearpage

\

\begin{figure}[h]
\centering
\includegraphics[trim=40.02bp 28.59bp 32.40bp 26.47bp,clip,width=\textwidth,height=\textheight,keepaspectratio]{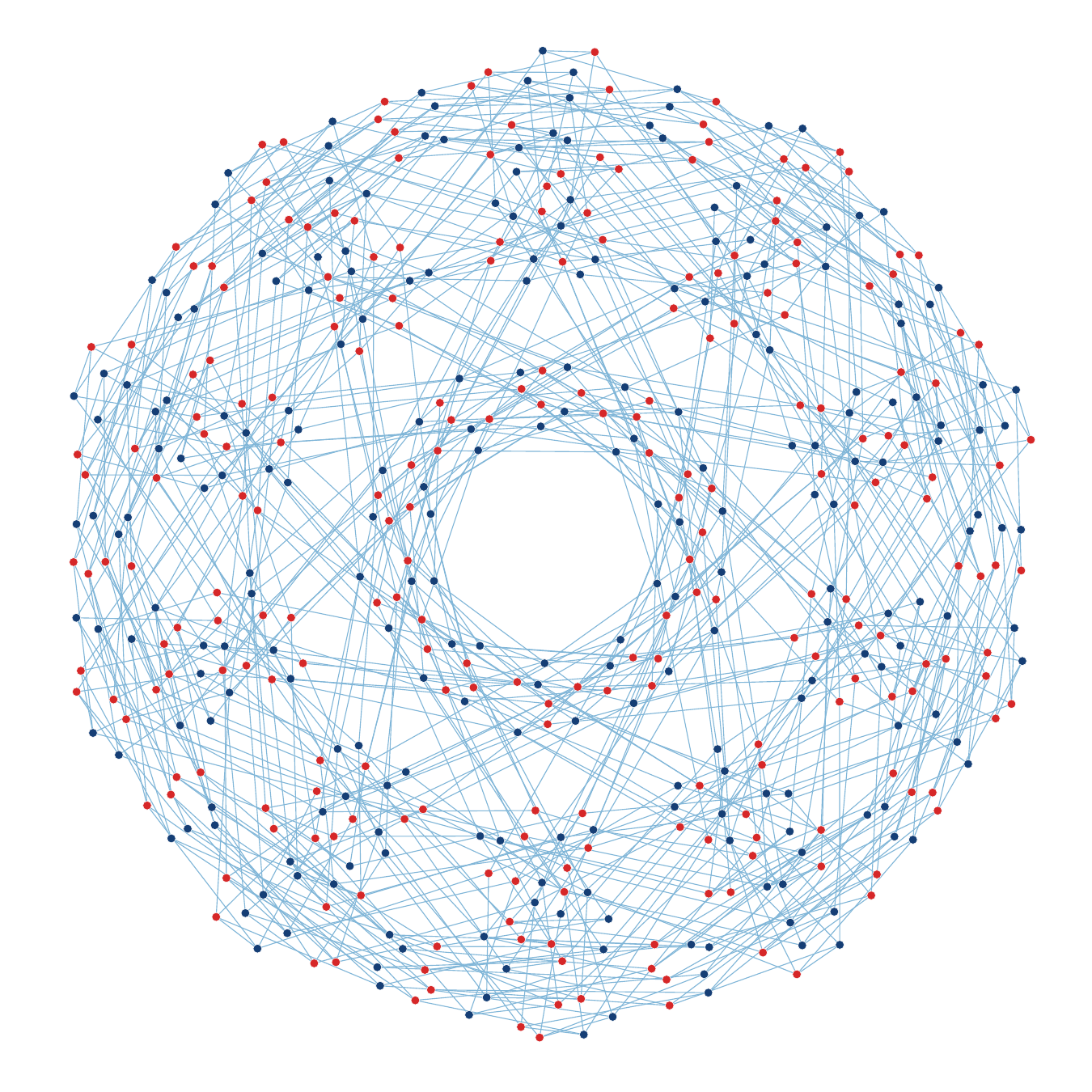}
\caption{An alternative drawing of the base graph $H$, with the same vertex colors.}
\end{figure}

\

\

\

\
\end{document}